\documentclass[11pt]{amsart}

\usepackage[ textheight=615pt,
  textwidth=360pt,
  centering
]{geometry}

\usepackage{xcolor}
\usepackage[english]{babel}
\usepackage[T1]{fontenc}
\usepackage[utf8]{inputenc}

\usepackage{amsmath,amssymb,amsthm,mathrsfs}
\usepackage{tikz-cd}
\usepackage{graphicx}
\usepackage{xypic}
\usepackage{hyperref}

\hypersetup{
  colorlinks=true,
  linkcolor=blue,
  urlcolor=red,
  pdftitle={}
}

\usepackage{graphicx}
\usepackage{xypic}
\usepackage{enumitem}
\usepackage{datetime}
\usepackage{xcolor}
\usepackage{fancyhdr}

\newtheorem{theorem}{Theorem}[section]
\newtheorem*{theorem*}{Theorem}

\newcommand{\X}{\mathcal{X}}

\newcommand{\m}{\omega}

\begin{document}

\title{On the coefficients of the asymptotic of analytic torsion}

\author{Mounir HAJLI}
\thanks{ }

\address{}

\email{hajlimounir@gmail.com}


            \date{\today, \currenttime}
            
            \maketitle 
            
\begin{abstract}
In this article, we study the asymptotic expansion of the analytic
torsion associated with positive Hermitian line bundles on projective
spaces. We show that the expansion contains no terms of the form
$k^{-i}\log k$ for $i\geq 1$. The proof combines the arithmetic
Riemann-Roch theorem with an explicit computation of the arithmetic
degree of spaces of global sections.
\end{abstract}

\tableofcontents

\section{Introduction}

Let $(M,g^{TM},\Theta)$ be a compact Hermitian manifold of complex
dimension $n$, and let $(L,h^L)$ be a positive Hermitian line bundle
over $M$. We denote by
\[
T(g^{TM},h^{kL})
\]
the analytic torsion associated with $g^{TM}$ and $h^{kL}$, where
$kL$ denotes $L^{\otimes k}$.

Finski \cite{Fin_JFA} established a full asymptotic expansion for the
analytic torsion $T(g^{TM},h^{kL})$. More precisely, there exist
coefficients $\alpha_i,\beta_i\in\mathbb R$, $i\in\mathbb N$, such that,
for every $\ell\in\mathbb N$,
\begin{equation}
\label{eq:Finski}
-T(g^{TM},h^{kL})
=
\sum_{i=0}^{\ell}
k^{n-i}
\left(
\alpha_i\log k+\beta_i
\right)
+
o(k^{n-\ell}),
\qquad k\to\infty.
\end{equation}
Moreover, the coefficients $\alpha_i$ are independent of the choices
of $g^{TM}$ and $h^L$.

In this article, we investigate the vanishing of the coefficients
$\alpha_i$ for positive Hermitian line bundles on projective spaces.
More precisely, our main result is the following.

\begin{theorem}\label{thm:vanishing}
Let $h^L$ be a smooth positive Hermitian metric on
$\mathcal O(1)$, and let $\mu$ be a smooth volume form on
$\mathbb P^n(\mathbb C)$. Then
\[
\alpha_i=0
\qquad\text{for all } i\geq n+1.
\]
\end{theorem}

\section{Analytic Torsion and the Arithmetic Riemann--Roch Theorem}

We briefly recall the definition of analytic torsion and the arithmetic
Riemann-Roch theorem of Gillet--Soulé. For details, we refer the reader
to \cite{Soule, Character, Character2, ARR, Ma}.

\medskip

\subsection*{Analytic torsion}

 Let $(M,g^{TM},\Theta)$ be a compact Hermitian manifold of complex
dimension $n$, and let $(L,h^L)$ be a smooth Hermitian line bundle over
$M$. Let $\zeta_q$ denote the zeta function associated with the Kodaira
Laplacian acting on smooth $(0,q)$-forms with values in $L^{\otimes k}$,
for $q=0,\ldots,n$; see \cite[Chapter VI, \S3]{Soule}. Following
\cite[p.~132]{Soule}, we define the analytic torsion of
$(M,g^{TM},\Theta)$ and $(L^{\otimes k},h^{kL})$ by
\[
T(g^{TM},h^{kL})
=
\sum_{q\geq 0}(-1)^{q+1}q\,\zeta_q'(0).
\]

Let $\lambda(kL)$ denote the determinant of cohomology of $kL$; see
\cite[Chapter VI]{Soule}. It is endowed with the Quillen metric
defined by
\[
h_Q=e^{T(g^{TM},h^{kL})}h_{L^2},
\]
where $h_{L^2}$ is the Hermitian metric on the determinant of
cohomology induced by the $L^2$-scalar product; see
\cite[p.~132]{Soule}. We denote the resulting Hermitian line bundle by
$(\lambda(kL),h_Q)$.

\medskip

We assume from now on that the metric $h^L$ is positive. Finski
\cite{Fin_JFA} established a full asymptotic expansion for
$T(g^{TM},h^{kL})$. More precisely, there exist coefficients
$\alpha_i,\beta_i\in\mathbb R$, $i\in\mathbb N$, such that, for every
$\ell\in\mathbb N$,
\begin{equation}
\label{eq:Finski}
-T(g^{TM},h^{kL})
=
\sum_{i=0}^{\ell}
k^{n-i}
\left(
\alpha_i\log k+\beta_i
\right)
+
o(k^{n-\ell}),
\qquad k\to\infty.
\end{equation}
Moreover, the coefficients $\alpha_i$ are independent of the choices
of $g^{TM}$ and $h^L$.

Finski's result generalizes the asymptotic result of
Bismut--Vasserot \cite[Theorem 8]{Bismut-Vasserot1}.
They proved that
\[
-T(g^{TM},h^{kL})
=
\alpha_0 k^n\log k+\beta_0 k^n+o(k^n),
\qquad k\to\infty,
\]
where
\[
\alpha_0
=
\frac{n}{2}
\int_M\frac{\omega^n}{n!},
\qquad
\beta_0
=
\frac12
\int_M
\log\left(
\det\frac{\widehat R^L}{2\pi}
\right)
\frac{\omega^n}{n!},
\]
and
\[
\omega=c_1(L,h^L).
\]

\medskip

When $\Theta=\omega$, Finski also computed the next coefficients:
\begin{equation}
\label{eq:alpha1}
\alpha_1
=
\frac{(3n+1){}}{12}
\int_M
c_1(TM)\frac{\omega^{n-1}}{(n-1)!}
\end{equation}
and
\begin{equation}
\label{eq:beta1}
\beta_1
=
\frac{{1}}{24}
\left(
24\zeta'(-1)+2\log(2\pi)+7
\right)
\int_M
c_1(TM)\frac{\omega^{n-1}}{(n-1)!}
\end{equation}

It is worth noting that the coefficient $\beta_1$ is independent of 
metrics.

\bigskip

\subsection*{The arithmetic Riemann-Roch theorem}

We now recall the arithmetic Riemann-Roch theorem. We use the standard
notation and refer the reader to \cite{AIT, Character, Character2} for
the foundations of arithmetic intersection theory. Let
\[
f:\mathcal X\longrightarrow\operatorname{Spec}(\mathbb Z)
\]
be the structural morphism of a smooth projective arithmetic variety
$\mathcal X$. Let $E$ be an algebraic vector bundle on $\mathcal X$.
We equip $E$ and $T\mathcal X$ with Hermitian metrics $h^E$ and
$h^{T\mathcal X}$, respectively, invariant under complex conjugation,
and assume that $h^{T\mathcal X}$ is induced by a Kähler form. By
\cite[Theorem 7]{ARR}, the first arithmetic Chern class of the
determinant of cohomology endowed with the Quillen metric
$(\lambda(E),h_Q)$ is
\begin{equation}
\label{eq:ARR-Pn}
\widehat c_1(\lambda(E),h_Q)
=
f_*
\left(
\widehat{\operatorname{ch}}(E,h^E)
\widehat{\operatorname{Td}}(T\X,h_{T\X})
-
a\left(
\operatorname{ch}(E_{\mathbb C})
\operatorname{Td}(T\X_{\mathbb C})
R(T\X_{\mathbb C})
\right)
\right)^{(1)}.
\end{equation}
The notation $\alpha^{(1)}$ denotes the component of degree one of
$\alpha\in\widehat{CH}^{\bullet}(\mathcal X)_{\mathbb Q}$, the arithmetic Chow ring, while
$a$ is the map from the real cohomology of $\mathcal X(\mathbb C)$ to
the arithmetic Chow ring defined in
\cite[3.3.4]{AIT} and \cite[2.2.1]{ARR}. The additive characteristic
class $R$ is associated with the power series
\[
R(x)
=
\sum_{\substack{m\geq1\\m\ {\rm odd}}}
\left(
2\zeta'(-m)
+
\left(1+\frac12+\cdots+\frac1m\right)\zeta(-m)
\right)
\frac{x^m}{m!}.
\]

\section{Proof of the Main Theorem \ref{thm:vanishing}}

Let $\mathbb P^n_{\mathbb Z}$ denote projective $n$-space over
$\operatorname{Spec}(\mathbb Z)$, regarded as an arithmetic variety.
We first consider the Fubini--Study metric and the standard arithmetic
model $\mathbb P^n_{\mathbb Z}$. For $k\in\mathbb N$, let $
E=\mathcal O(k)$
and equip it with the Fubini-Study metric. We denote the resulting
Hermitian line bundle by $\overline{\mathcal O(k)}_{\mathrm{FS}}$. This data
defines a metrized $\mathbb Z$-module structure on
\[
H^0(\mathbb P^n_{\mathbb Z},\mathcal O(k)).
\]
The arithmetic degree of this Hermitian vector bundle $\overline{H^0\left(\mathbb P^n_{\mathbb Z},\mathcal O(k)\right)}_{\mathrm{FS}}$ is given by

\begin{equation}
\label{eq:asympFS}
\widehat{\deg}
\left(
\overline{H^0\left(\mathbb P^n_{\mathbb Z},\mathcal O(k)\right)}_{\mathrm{FS}}
\right)
=
-\frac12
\sum_{m_0+\cdots+m_n=k}
\log
\left(
\frac{m_0!\cdots m_n!\,n!}{(k+n)!}
\right).
\end{equation}
This formula follows from the explicit description of the $L^2$-metric
on the space of homogeneous polynomials; see, for example,
\cite[\S1.3]{Ran1}.

With the notation of \cite{Ran1}, set
\[
\eta_n^{(k)}
:=
-\frac12
\sum_{m_0+\cdots+m_n=k}
\log
\left(
\frac{m_0!\cdots m_n!\,n!}{(k+n)!}
\right).
\]
The Euler--Maclaurin expansion obtained in
\cite[Proposition 4.2.3]{Ran1} implies that, for every
integer $r>0$, $\eta_n^{(k)}$ admits an asymptotic expansion
\begin{equation}
\label{eq:eta-expansion}
\begin{aligned}
\eta_n^{(k)}
={}&
a_{n+1}k^{n+1}
+b_nk^n\log k
+a_nk^n
+b_{n-1}k^{n-1}\log k
+\cdots
\\
&\qquad
+b_0\log k+a_0
+c_{-1}k^{-1}
+c_{-2}k^{-2}
+\cdots
+c_{-r}k^{-r}
+O(k^{-r-1}),
\end{aligned}
\end{equation}
as $k\to\infty$.

The first coefficients are
\[
a_{n+1}
=
\frac{\sigma_n}{(n+1)!},
\qquad
b_n
=
\frac{1}{4(n-1)!},
\]
and
\[
a_n
=
\frac{1}{2n!}
\left(
(n+1)\sigma_n
-\log(n!)
-\frac n2\log(2\pi)
\right),
\]
where $\sigma_n$ denotes the Stoll number. More generally, for
$0\leq i\leq n$,
\[
a_i
\in
\mathbb Q
+
\mathbb Q\log(n!)
+
\mathbb Q\zeta'(0)
+\cdots+
\mathbb Q\zeta'(-n+i),
\]
while
\[
b_i\in\mathbb Q,
\]
and, for every $i>0$,
\[
c_{-i}\in\mathbb Q.
\]

Applying the arithmetic Riemann--Roch theorem \eqref{eq:ARR-Pn} to
$E=\overline{\mathcal O(k)}_{\mathrm{FS}}$ and using the multiplicativity of the
arithmetic Chern character, we obtain the arithmetic Riemann--Roch
formula for $\overline{\mathcal O(k)}_{\mathrm{FS}}$.

For $k$ sufficiently large, the ampleness of $\mathcal O(1)$ implies
the vanishing of the higher cohomology groups. Hence, by the definition
of the Quillen metric,

\begin{equation}
\label{eq:Quillen-Pn}
\widehat{\deg}
\left(
\lambda({\overline{\mathcal O(k)}_{\mathrm{FS}}}),h_Q
\right)
=
\widehat{\deg}
\left(
H^0(\mathbb P^n_{\mathbb Z},\mathcal O(k))_{\mathrm{FS}}
\right)
-
\frac12
T(g^{T\mathbb P^n}_{\mathrm{FS}},{\overline{\mathcal O(k)}_{\mathrm{FS}}}).
\end{equation}

The right-hand side of the arithmetic Riemann--Roch formula is a
polynomial in $k$. Consequently,
\[
\widehat{\deg}\bigl(\lambda(\mathcal O(k)_{\mathrm{FS}}),h_Q\bigr)
\]
is a polynomial in $k$.

Combining this observation with \eqref{eq:Quillen-Pn} and the expansion
\eqref{eq:eta-expansion}, we conclude that the asymptotic expansion of
\[
T(g^{T\mathbb P^n}_{\mathrm{FS}},
  \overline{\mathcal O(k)}_{\mathrm{FS}})
\]
can not contain any term of the form
\[
k^{-j}\log k,\qquad j\geq1.
\]

On the other hand, Finski's expansion  \eqref{eq:Finski}  contains such terms
precisely through the coefficients
\[
\alpha_{n+j}k^{-j}\log k.
\]
Consequently,
\[
\alpha_{n+j}=0
\qquad\text{for every }j\geq1.
\]
This proves
\[
\alpha_i=0
\qquad\text{for all }i\geq n+1.
\]

Finally, the coefficients $\alpha_i$ are independent of the choice of
Hermitian metrics by Finski's theorem. Hence the conclusion holds for
every smooth positive Hermitian metric on $\mathcal O(1)$.\\

\bigskip

We are currently investigating whether the same argument extends to
polarized toric varieties equipped with canonical metrics in the sense
of Guillemin \cite{Guillemin1}; see also \cite{Abreu1}. In this setting,
the relevant arithmetic degrees admit explicit combinatorial
descriptions.

\section*{Acknowledgements}

The author would like to thank Professor Siarhei Finski for his helpful comments
on an earlier version of this article.

\bibliographystyle{plain} 

\bibliography{biblio}

\end{document}